\documentclass{article}

\usepackage{amsmath,amssymb,latexsym,psfrag}
\usepackage{amsfonts}
\usepackage{graphicx}
\usepackage{epsf,epsfig,color}

\newtheorem{theo}{Theorem}

\newtheorem{lemm}{Lemma}
\newtheorem{defn}{Definition}

\def\K{\mathbb{K}}

\def\E{\mathbb{E}}
\def\S{\mathbb{S}}
\def\0{{\bf 0}}

\def\R{\mathbb{R}}

\def\B{{ B}}

\def\w{{w'}}
\def\wo{{{w_0}'}}

\def\en2{{\epsilon_n^2}}
\renewcommand{\E}{\mathbb E \,}
\newcommand{\C}{{\cal C}}

\newcommand{\tod}{\stackrel{{\cal D}}{\longrightarrow}}
\newcommand{\eqd}{\stackrel{{\cal D}}{=}}

\newcommand{\eqco}{\setcounter{equation}{0}}
\newcommand{\thco}{\setcounter{theo}{0}}
\newcommand{\prco}{\setcounter{prop}{0}}
\newcommand{\laco}{\setcounter{lemm}{0}}
\newcommand{\coco}{\setcounter{coro}{0}}
\newcommand{\cjco}{\setcounter{conj}{0}}

\newcommand{\deco}{\setcounter{defn}{0}}
\newcommand{\allco}{\eqco  \thco \prco \laco \coco \cjco \deco}
\newcommand{\X}{{\cal X}}
\def\b{{\beta}}

\def\e{{\epsilon_\la}}
\def\2e{{\epsilon^2_\la}}
\def\la{{\lambda}}

\def\v{{ w }}

\newcommand{\Y}{{\cal Y}}

\renewcommand{\P}{{{\cal P}}}

\newcommand{\Cov}{{\rm Cov}}

\newcommand{\Var}{{\rm Var}}

\newcommand{\Vertices}{{\rm Vertices}}

\def\bdm{\begin{displaymath}}
\newcommand{\edm}{\end{displaymath}}
\def\benu{\begin{enumerate}}
\def\eenu{\end{enumerate}}
\def\beqn{\begin{equation}}
\def\eeqn{\end{equation}}
\def\be{\begin{equation}}
\def\ee{\end{equation}}
\def\bea{\begin{eqnarray}}
\def\eea{\end{eqnarray}}
\newcommand{\bean}{\begin{eqnarray*}}
\newcommand{\eean}{\end{eqnarray*}}
\newcommand{\bear}{\begin{eqnarray}}
\newcommand{\eear}{\end{eqnarray}}

\def\Comment#1{
\marginpar{$\bullet$\quad{\tiny #1}}}

\def\R{\mathbb{R}}
\def\B^2{\mathbb{D}}
\def\S{\mathbb{S}}
\def\B{\mathbb{B}}

\def\b{{\beta}}

\def\th{{\theta}}

\def\txi{{\xi^{\la, z} }}

\def\de{{\delta}}

\def\T{{T^{\la, z}}}
\def\K{{K^{\la, z}}}
\def\Pl{{ {\cal P}^{\la, z}}}
\def\Plz{{ {\cal P}^{\la, z}_{r_z} }}
\def\pK{{ \partial K }}
\def\c{{ c^{\la, z} }}
\def\y{ \bar{y}}

\def\Comment#1{\lineskip-4pt
\marginpar{ $\bullet$\quad{\em\small #1}}}

\def\qed{\hfill\hbox{${\vcenter{\vbox{
    \hrule height 0.4pt\hbox{\vrule width 0.4pt height 6pt
    \kern5pt\vrule width 0.4pt}\hrule height 0.4pt}}}$}}

\def\la{{\lambda}}
\def\ka{{\kappa}}

\begin{document}
\title{\bf Variance asymptotics for random polytopes in smooth convex bodies}

\author{Pierre Calka$^{*}$,  J. E. Yukich$^{**}$ }

\maketitle

\footnotetext{{\em American Mathematical Society 2000 subject
classifications.} Primary 60F05, Secondary 60D05} 

\footnotetext{$~^{*}$ Research partially supported by French ANR grant PRESAGE.}\footnotetext{$~^{**}$ Research supported in part by NSF grant
DMS-1106619}

\begin{abstract} Let $K \subset \R^d$ be a smooth convex set and let $\P_\la$
be a Poisson point process on $\R^d$ of intensity $\la$.  The convex
hull of $\P_\la \cap K$ is a random convex polytope $K_\la$.  As
$\la \to \infty$, we show that the variance of the number of
$k$-dimensional faces of $K_\la$, when properly scaled, converges to
a scalar multiple of the affine surface area of $K$. Similar
asymptotics hold for the variance of the number of $k$-dimensional
faces
for the convex hull of a binomial process in $K$.
\end{abstract}

\section{Introduction}\label{INTRO}


 Let $K \subset \R^d$ be a compact convex body with non-empty interior and having a $C^3$ boundary
of positive Gaussian curvature $\ka$.   Letting $\P_\la$ be a
Poisson point process in $\R^d$ of intensity $\la$ we denote by
$K_\la$  the convex hull of $\P_\la \cap K$. Let $f_k(K_\la), k \in
\{0,1,...,d-1\},$ be the number of $k$ faces  of $K_\la$.

 R\'enyi and Sulanke \cite{RS} were the first to consider the
average behavior of $f_0(K_\la)$ in the planar case. Generalizing their formula to higher dimensions, B\'ar\'any \cite{Ba} showed there is a
constant $D_{0,d}$ such that
$$\lim_{\la \to \infty} \la^{-(d-1)/(d + 1)} \E f_0(K_\la) = D_{0,d}
\int_{\partial K} \ka(z)^{1/(d + 1)} dz.
$$
The integral $\int_{\partial K} \ka(z)^{1/(d + 1)} dz$ is 
 known as the {\em affine surface area}
of $\pK$.
Assuming only that $K$ has a boundary $\partial K$  of
differentiability class $C^2$, Reitzner \cite{Re1} extended this
result to $f_k(K_\la), \ k \in \{0,1,...,d-1\}$, showing for all $d
\geq 2$  that there are constants $D_{k,d}$ such that \be
\label{ReLLN} \lim_{\la \to \infty} \la^{-(d-1)/(d + 1)} \E
f_k(K_\la) = D_{k,d} \int_{\partial K} \ka(z)^{1/(d + 1)} dz. \ee

Reitzner \cite{Re} also showed that $(f_k(K_\la)- \E
f_k(K_\la))/\sqrt{ \Var f_k(K_\la) }$  converges in distribution to
a mean zero normal random variable as $\la \to \infty$, though there
have been relatively few results concerning the asymptotic variance
of $f_k(K_\la)$. Theorem 4 of Reitzner \cite{Re} gives upper and
lower bounds of the same magnitude for $\Var f_k(K_\la), k \in
\{0,1,...,d - 1\},$ which extends work of Buchta \cite{Bu}, who
obtains lower bounds for $\Var f_0(K_\la)$ of order $\la^{(d-1)/(d +
1)}$. In the special case that $K$ is a ball, closed form variance
asymptotics for $\Var f_k(K_\la), k \in \{0,1,...,d - 1\}$ are given
in \cite{SY, CSY}.

Let  $K_n'$ be the convex hull of $n$ i.i.d. random variables
uniformly distributed on $K$. Our main two results 
resolve the open question of determining
variance asymptotics for $\Var f_k(K_\la)$ and $\Var f_k(K_n')$, $K$
smooth and convex,  as put forth on
page 1431 of  \cite{WW}. 

\begin{theo} \label{Th1} For all $k \in \{0,1,...,d-1\}$, there
exist positive constants $F_{k,d}$  such that  \be \label{main1}
\lim_{\la \to \infty} \la^{-(d-1)/(d + 1)} \Var f_k(K_\la) = F_{k,d}
\int_{\partial K} \ka(z)^{1/(d + 1)} dz. \ee
\end{theo}

Let $\rm{vol}$ be the Lebesgue measure. De-Poissonization methods,
based on coupling, yield  the following binomial counterpart of
\eqref{main1}. When $k = 0$, it  resolves Conjecture 1 of Buchta
\cite{Bu}.
\begin{theo} \label{Th1b} For all $k \in \{0,1,...,d-1\}$ we have \be \label{main1b}
\lim_{n \to \infty} n^{-(d-1)/(d + 1)} \Var f_k(K_n') = F_{k,d}({\rm{vol}}(K))^{-(d-1)/(d+1)}
\int_{\partial K} \ka(z)^{1/(d + 1)} dz. \ee
\end{theo}

\noindent{\em Remarks.}

(i) {\em Related work.} B\'ar\'any and Reitzner (page 3 of \cite{BR})
conjecture for general convex bodies that $\Var f_k(K_\la)$ should -
up to constants  - behave like the variance of the volume of the wet
part of the floating body, which, in the case of smooth convex sets,
is proportional to affine surface area.  Theorem \ref{Th1} resolves
a sharpened version of this conjecture in the case that $\partial K$
is smooth.

(ii) {\em The constants $F_{k,d}$}. The proofs of Theorems \ref{Th1}
and \ref{Th1b} show that $F_{k,d}$ is defined in terms of parabolic
growth processes on the upper half-space $\R^{d-1} \times \R^+$.  As
noted on page 137 of Buchta \cite{Bu}, $F_{k,d}$ may also be
identified in terms of a constant involving complicated double
integrals given in Groeneboom \cite{Gr}.

(iii)  {\em Volume asymptotics}.  Under a $C^3$ and $C^2$ assumption
on $\pK$, respectively, B\'ar\'any \cite{Ba} and Reitzner \cite{Re2}
show \be \label{emainint} \lim_{\la \to \infty} \la^{2/(d + 1)} \E
{\rm{vol}}(K \setminus K_\la) = c_{d} ({\rm{vol}}(K))^{2/(d + 1)}
\int_{\partial K} \ka(z)^{1/(d + 1)} dz. \ee B\"or\"oczky et al.
\cite{BFH} extend this limit and \eqref{ReLLN} to convex hulls  of
i.i.d. points having a non-uniform density on $K$. Theorem 3 of
Reitzner's breakthrough paper \cite{Re} gives upper and lower bounds
of the same magnitude for $\Var\; {\rm{vol}}(K_\la)$, though it falls
short of giving 
 a limiting
variance. Notice that Theorems \ref{Th1} and \ref{Th1b} fill in this
gap as follows.
Buchta notes (see Corollary 1 and (3.6) of \cite{Bu}) under
sufficient smoothness of $\partial K$,  that variance asymptotics
for $n^2 \Var f_0(K_n')$ coincide with variance asymptotics for
$\Var\; {\mathrm{ vol}}(K_n')$, that is
$$\Var\; {\rm{vol}}(K_n')=\frac{\Var(f_0(K_{n + 2}'))+d_{n+2}}{(n+1)(n+2)},$$
where $$\lim_{n\to\infty}\left(\frac{3-d}{d+1}\int_{\partial
K}\kappa(z)^{1/(d + 1)}dz\cdot n^{(d-1)/(d+1)}\right)^{-1}d_n=1.$$
Consequently, putting $G_d:=F_{0,d}+(3-d)/(d+1)$ and putting $k = 0$
in \eqref{main1b}, we get
  \be \label{mainvol}
\lim_{n \to \infty} n^{(d+3)/(d + 1)} \Var\; {\mathrm{vol}}(K_n') =
G_{d}({\rm{vol}}(K))^{(d+3)/(d+1)} \int_{\partial K} \ka(z)^{1/(d +
1)} dz.\ee By \eqref{mainvol} and Proposition 3.2 of \cite{Vu},
which states that $\Var\;{\rm{vol}}(K_n')$ and $\Var\;{\rm{vol}}(K_n)$
coincide up to first order, 
we deduce
 \be \label{mainvolb}
\lim_{\lambda \to \infty} \lambda^{(d+3)/(d + 1)} \Var\; {\mathrm{vol}}(K_\la) = G_{d}
\int_{\partial K} \ka(z)^{1/(d + 1)} dz. \ee


\ \




The paper is organized  as follows. 
Section 2 introduces the main tool for the proof of Theorem
\ref{Th1}, 
 namely  the paraboloid
  growth process used  in \cite{SY} and \cite{CSY}.
We state a general result, Theorem \ref{Th2},  giving expectation
 and variance asymptotics for the empirical
$k$-face measure, which includes Theorem \ref{Th1} as a special
case. Theorem \ref{Th2} also shows that the constants $F_{k,d}$ of
Theorem \ref{Th1} may be expressed in terms of integrals of one and
two point correlation functions of a scaling limit $k$-face
functional $\xi_k^{(\infty)}$ associated with parabolic growth
processes. Section 3 introduces {\em an affine transform of $K$ and
a scaling transform of the affine transform} to link  the finite
volume $k$ face functional with its infinite volume scaling limit
counterpart $\xi_k^{(\infty)}$.
 Section 4 contains the main technical aspects  of the paper, focussing on the properties
 of the  re-scaled $k$-face functionals. 
In particular  Lemmas  \ref {L4.2} and \ref{L4.5} show that the one
and two point correlation functions of the re-scaled $k$-face
functional on the affine transform of $K$ are well approximated by
the corresponding one and two point correlation functions of the
re-scaled $k$-face functional on an osculating ball. In this way the
expectation and variance asymptotics for $f_k(K_\la)$, $K$ an
arbitrary smooth body,  are controlled by the corresponding
asymptotics for $f_k(K_\la)$ when $K$ is a ball. The latter
asymptotics are established in \cite{CSY}.
  Section 5 contains the proof of Theorem \ref{Th2} which implies Theorem \ref{Th1}.
  Finally, in Section 6, we prove the  de-Poissonized limit  \eqref{main1b}. 



\section{Paraboloid growth processes and a general result}

\allco

Given a finite point set $\X \subset \R^d$, let
$\rm{co}(\X)$ be its convex hull.

\begin{defn} \label{kface} Given $k \in \{0,1,...,d-1\}$
and $x$ a vertex of $\rm{co}(\X)$, define the $k$-face
functional $\xi_k(x, \X)$ to be the product of $(k +1)^{-1}$ and
the number of $k$ faces of $\rm{co}(\X)$ which contain $x$.
Otherwise we put $\xi_k(x, \X) = 0$.
The empirical k-face measure is \be \label{zerom}
\mu^{\xi_k}_\la:= \sum_{x \in \P_\la \cap K} \xi_k(x, \P_\la
\cap K) \de_x, \ee where $\de_x$ is the unit point mass at $x$.
\end{defn}

Thus the number of $k$-faces in $\rm{co}(\X)$ is $\sum_{x \in \X}
\xi_k(x, \X)$. We shall give a general result describing the limit
behavior of $\mu^{\xi_k}_\la$ in terms of parabolic growth processes on
$\R^d$.

\vskip.5cm

\noindent{\bf Paraboloid growth processes}. Denote points in $\R^{d-1} \times
\R$ by $\v:= (v,h)$ or $\v':= (v',h')$, depending on context. Let
$\Pi^{\uparrow}$ be the epigraph of the parabola $v \mapsto
|v|^2/2$, that is $\Pi^{\uparrow} := \{(v,h) \in \R^{d-1} \times
\R^+, h \geq |v|^2/2\}.$  Letting $\X \subset \R^d$ be locally
finite,   define the parabolic growth model
$$
\Psi(\X):= \bigcup_{ \v \in \X } (\v \oplus \Pi^{\uparrow}),
$$
where $\oplus$ denotes Minkowski addition. A point $w_0 \in \X$ is
{\em extreme} with respect to $\Psi(\X)$ if the epigraph $w_0 \oplus
\Pi^{\uparrow}$ is not a subset of the union of the epigraphs $\{w
\oplus \Pi^{\uparrow}, w \in \X \setminus w_0\}$, that is $(w_0
\oplus \Pi^{\uparrow}) \nsubseteq \bigcup_{\v \in \X \setminus w_0}
(w \oplus \Pi^{\uparrow}).$

The paraboloid hull model $
\Phi(\X)$ is defined  as in Definition 3.4 of \cite{CSY}:
$$
\Phi(\X):=\bigcup_{\left\{\substack{w\in \R^{d-1}\times\R\\(w\oplus\Pi^{\downarrow})\cap \X=\emptyset}\right.}
(w\oplus\Pi^{\downarrow}),
$$
where $\Pi^{\downarrow}:= \{(v,h) \in \R^{d-1} \times \R, h \le-
|v|^2/2\}.$
 It may be viewed as the dual of the paraboloid growth model
$\Psi(\X)$. Let $\P$ be a rate one homogeneous Poisson point process
on $\R^{d-1} \times \R^+$ and let  $\Psi:= \Psi(\P)$  and $\Phi:=
\Phi(\P)$ be the corresponding {\em paraboloid growth and hull
processes}.
 As in \cite{CSY}, the set $\Vertices(\Phi)$ coincides with the extreme points
of $\Psi$.

\begin{defn} \label{xiinf} (cf. section 6 of \cite{CSY}). Define the
scaling limit $k$-face functional $\xi^{(\infty)}_{k}(x,
\P)$,  for $x \in \P$, and $k \in \{0, 1,..., d - 1\}$, to be the
product of $(k + 1)^{-1}$ and the number of $k$-dimensional
paraboloid faces of the hull process $\Phi$ which contain $x$,
if $x$ belongs to $\Vertices(\Phi)$, and zero otherwise.
\end{defn}

One of the main features of our approach is that
$\xi^{(\infty)}_{k}$ is indeed a scaling limit of appropriately
re-scaled $k$-face functionals, as seen in Lemma \ref{L4.3} of
Section 4 and also in Lemma 7.2 of \cite{CSY}.

Define the following second order correlation functions for
$\xi^{(\infty)}(x, \P) :=\xi^{(\infty)}_{k}(x, \P)$ (cf. (7.2),
(7.3) of \cite{CSY}).

\begin{defn} For all $\v_1, \v_2 \in \R^d$, put
\be \label{SO1} \zeta_{\xi^{(\infty)}}(\v_1):=
\zeta_{\xi^{(\infty)}}(\v_1, \P):= \E \xi^{(\infty)}(\v_1, \P)^2 \ee
and \be \label{SO2} \zeta_{\xi^{(\infty)}}(\v_1,\v_2):=
\zeta_{\xi^{(\infty)}}(\v_1, \v_2, \P):= \ee $$
 \E \xi^{(\infty)}(\v_1, \P
\cup \{\v_2\} ) \xi^{(\infty)}(\v_2, \P \cup \{\v_1\} ) -  \E
\xi^{(\infty)}(\v_1, \P ) \E \xi^{(\infty)}(\v_2, \P ). $$
\end{defn}
Note that \be \label{S03}
 \sigma^2(\xi^{(\infty)}) := \int_0^{\infty} \varsigma_{\xi^{(\infty)}}((\0,h)) dh +
   \int_0^{\infty}  \int_{\R^{d-1}} \int_0^{\infty}
   \varsigma_{\xi^{(\infty)}}((\0,h),(v',h')) dh' dv' dh
\ee is finite and positive by Theorems 7.1 and 7.3 in \cite{CSY}.

\vskip.5cm

Theorem \ref{Th1} is a special case of the following more general
result giving the asymptotic behavior of the empirical k-face
measures in terms of parabolic growth processes.  Let $\C(K)$ be the
class of continuous functions on $K$ and let $\langle g, \mu_\la^\xi
\rangle$ denote the integral of $g$ with respect to $\mu_\la^\xi$.

\begin{theo} \label{Th2}  For all $g \in \C(K)$ and $k \in \{0, 1,..., d - 1\}$, we have
\be \label{main3} \lim_{\la \to \infty} \la^{-(d-1)/(d + 1)} \E
[\langle g, \mu_\la^{\xi_k} \rangle] = \int_0^\infty \E
\xi^{(\infty)}_{k}((\0,h), \P) dh \int_{\pK} g(z) \ka(z)^{1/(d +
1)} dz \ee and \be \label{main4} \lim_{\la \to \infty}
\la^{-(d-1)/(d + 1)} \Var [\langle g, \mu_\la^{\xi_k} \rangle] =
\sigma^2(\xi^{(\infty)}_{k}) \int_{\pK} g(z)^2 \ka(z)^{1/(d + 1)}
dz. \ee
\end{theo}


{\em Remarks}.

(i) {\em Related work}. Up to now, \eqref{main4} has
been known only for bodies of constant curvature, i.e., only for $K
= r\B^d, d \geq 2, r > 0$; see Theorem 7.3 of \cite{CSY}.

(ii) {\em The constants}. Recalling the notation of Theorem
\ref{Th1}, we obtain $F_{k,d} = \sigma^2(\xi^{(\infty)}_{k})$.

(iii) {\em Extensions}. As in  \cite{Re} and \cite{BFH}, we expect
that the  $C^3$ boundary condition could be relaxed to a $C^2$
condition, and we comment on this in Section 5.3. Following the
methods of Section 6, we obtain the counterpart of Theorem \ref{Th2}
for binomial input.

\section{Affine and scaling transformations} 

\allco

For each $z \in \partial K$, we first consider
an affine transformation  ${\cal A}_z$  of  $K$, one under which the scores $\xi_k$
are invariant, but under which the principal curvatures of ${\cal
A}_z(K)$ at $z$ coincide, that is to say ${\cal
A}_z(K)$ is `umbilic' at $z$.  This property allows us to readily approximate the
functionals $\xi_k$ on Poisson points in ${\cal A}_z(K)$ by the
corresponding functionals on Poisson points in the `osculating ball'
at $z$, defined below. The key idea of replacing the mother body $K$ with an osculating ball
 has been used by R\'enyi and Sulanke \cite{RS2}, B\'ar\'any
 \cite{Ba}, and B\"or\"oczky et al.
\cite{BFH}, among others.

We in turn transform ${\cal A}_z(K)$ to a subset of $\R^{d-1} \times
\R$ via scaling transforms $T^{\la, z}, \la \geq 1$.  These
transforms yield re-scaled $k$-face functionals $\xi^{\la, z}$ on
the Poisson points $T^{\la, z}( \P_\la \cap {\cal A}_z(K))$, ones
 which are well approximated
by re-scaled $k$-face functionals on the image under $T^{\la, z}$ of Poisson points
in the osculating ball
at $z$.  In the large $\la$ limit the latter in turn converge to
the scaling
limit functionals $\xi^{(\infty)}$ given in Definition \ref{xiinf}.

In this way  the expectation and variance
asymptotics for $k$-face functionals on Poisson points in $K$ are
obtained by averaging, with respect to all $z \in
\partial K$, the respective asymptotics for the re-scaled $k$-face functionals on
Poisson points in osculating balls at $z$.  The limit theory of the latter is established
in \cite{CSY, SY} and we shall draw upon it in our approach.


\vskip.5cm

\noindent{\bf 3.1. Affine transformations ${\cal A}_z, z \in K$.}
Let ${\mathcal M}(K)$ be the medial axis of $K$.  ${\mathcal M}(K)$
has Lebesgue measure zero and we parameterize points $x \in K
\setminus {\mathcal M}(K)$ by $x := (z, t)$, where $z \in
\partial K$ is the unique boundary point closest to $x$ and where $t
\in [0, \infty)$ is the distance between
$x$ and $z$.  


 Denote by $C_{z,1},\cdots,C_{z,d-1}$ the
principal curvatures of $\partial K$ at $z$, i.e. the eigenvalues of
the Weingarten operator at $z$. Let $\ka(z): = \prod_{i=1}^{d - 1}
C_{z,i}$ be the Gaussian curvature at $z$, so that the {\em Gaussian
curvature radius $r_z$}
satisfies $\kappa(z)= r_z^{-(d-1)}$. 

For $z \in \partial K$, consider the  affine transformation
${\mathcal A}_z$ which preserves $z$, the Lebesgue measure, the unit
inner normal to $z$, and which transforms the Weingarten operator at
$z$ into $r_z^{-1}I_{d-1}$ where $I_{d-1}$ is the identity matrix of
$\R^{d-1}$. Under the action of ${\mathcal A}_z$, the number of
$k$-faces of the random convex hull inside the mother body $K$ is
preserved.
 Additionally,  $\xi_k, k \in \{0,1,...,d-1\}$ is stable under
the action of ${\mathcal A}_z$, namely \be \label{stable} \xi_k(x,
\P_\la \cap K) = \xi_k({\mathcal A}_z(x), {\mathcal A}_z(\P_\la \cap
K)). \ee  Indeed, ${\mathcal A}_z$ sends any $k$-face of $K_{\la}$
to a $k$-face of ${\mathcal A}_z(K_{\la})$. This follows since
affine transformations
 preserve convexity and convex hulls. A $k$-face $F_k
$ of $K_{\la}$ is a.s. the convex hull of $(k+1)$ points from $\P_\la$, so it is sent to the
 convex hull of the images by ${\mathcal A}_z$. Moreover, any support hyperplane $H$ such that $H\cap K_{\la}=F_k$
 is sent to a support hyperplane of the image of $K_{\la}$ such that its intersection with it is the image of the
 face $F_k$. So the image of $F_k$ is also a $k$-face of the image of $K_{\la}$.

 Put
$K_z:= {\mathcal A}_z (K)$. By construction the principal curvatures at $z$
all equal $r_z^{-1}$. We recall that ${\mathcal A}_z$
preserves the distribution of ${\mathcal P}_\la$ so in the sequel,
we will make a small abuse of notation by identifying ${\mathcal
P}_{\la}$ and $K_{\la}$ with ${\mathcal A}_z({\mathcal P}_{\la})$
and ${\mathcal A}_z(K_{\la})$, respectively.
Define the {\em osculating ball} at $z \in \partial K$ to be the
ball whose center $z_0:=z_0(z)$ is at distance $r_z$ from $z$ along the
inner normal to $z$. Lemma \ref{L3.2} shows that the boundary of the
osculating ball $B_{r_z}(z_0)$ is not far from $\partial K_z$,
justifying the terminology.

Given $z \in
\partial K$, define $f:{\mathbb
S}^{d-1}\mapsto \R^+$ to be the function such that for all $u\in
{\mathbb S}^{d-1}$, $(z_0+f(u)u)$ is the point of the half line
$(z_0+\R^+u)$ contained in $\partial K_z$ and furthest from $z_0$.
Thus $\partial K_z$ is given by $(
f(u),u), u \in {\mathbb S}^{d-1}$. Given $z \in
\partial K$ we let the inner unit normal be $k_z:={ (z_0-z) /
|z-z_0| }$. Here and elsewhere we let $|w|$ denote the Euclidean
norm of $w$. For each fixed $z \in \pK$,
 we  parameterize points $\v$ in $\R^{d-1} \times \R$ by $(r, u)$ where $r := |\v - z_0|$ and where
  $u\in {\mathbb S}^{d-1}$. Henceforth, points $(r,u)$ are with
  reference to $z$.
 For $z=(r_z,u_z) \in \pK$ let $T_z \sim \R^{d-1}$ denote the tangent space to
 ${\mathbb S}^{d-1}$ at $u_z$. The exponential map on the sphere $\exp_{d-1}: T_z \to  {\mathbb S}^{d-1} $ maps a vector
 $v$ of the tangent space to the point $u \in  {\mathbb S}^{d-1} $ such that $u$ lies at the end of
 the geodesic of length $|v|$ starting at $z$ and having direction $v$.  We let the origin of the tangent space be at
 $u_z$.

 \vskip.5cm

\noindent{\bf 3.2. Scaling transformations $\T, z \in \partial K, \la
\geq 1$}.   Having transformed $K$ to $K_z$,  we now re-scale $K_z$ for all $\la \geq 1$ with
a scaling transform denoted $\T$. Our choice of $\T$ is motivated by the
following desiderata. First, consider the epigraph of $s_\la: \S^{d-1} \mapsto \R$
defined by
$$
s_\la(u, \P_\la) = r_z - h_{K_\la}(u), \ u \in \S^{d-1},
$$
where we recall that $r_z$ is the Gaussian curvature radius at $z$ and $h_{K_\la}(u):= \sup\{ \langle x, u \rangle, \ x \in K_\la\}$ denotes the support function of $K_\la$.
Noting  that $h_{K_\la}(u) = \sup_{x \in \P_\la} h_x(u)$ for all $u \in \S^{d-1}$, it follows that the considered
 epigraph is the union of epigraphs, which, locally near the vertices of $K_\la$, are of parabolic structure.
Thus any scaling transform should preserve this structure, as should the scaling limit.  Second,
a subset of $K_z$ close to $z$ and having a unit volume scaling image should host on average
$\Theta(1)$ points of $\P^{\la, z}$, that is to say the intensity density of the re-scaled points should
be of order  $\Theta(1)$.  As in Section 2 of \cite{CSY}, it follows that the transform $\T$ should re-scale $K_z$
in the $(d-1)$ tangential directions with factor $\la^{1/(d + 1)}$ and in the radial direction with
factor $\la^{2/(d + 1)}$.  It is easily checked that the following choice of $\T$ meet these
criteria; cf.  Lemma \ref{L3.3} below.  Throughout we put $$\b := {1\over d + 1}.$$

Define for all  $z \in \pK$ and $\la \geq 1$ the finite-size
scaling transformation $\T: \R^+ \times {\mathbb S}^{d -1} \to \R^{d -1}
\times \R$
by \be \label{Tdef} \T((r, u)) := \left(  (r_z^d\la)^\b
\exp_{d-1}^{-1}(u), (r_z^d\la)^{2\b}(1 - {r \over r_z})\right):= (v',h'):=
\w. \ee  Here $\exp_{d-1}^{-1}(\cdot)$ is the inverse exponential map, which is well defined
 on $\S^{d-1} \setminus \{ - u_z \}$ and which takes values in the ball of radius $\pi$ and centered at the origin of the  tangent space $T_z$. 
 We shall write
$v':=(r_z^d\la)^\b \exp_{d-1}^{-1}u:= (r_z^d\la)^\b v,$
where $v \in \R^{d-1}.$  
We put
$$
\T(K_z):= \K; \ \ \T(B_{r_z}(z_0)):= B^{\la, z};  $$ $$  \T(\P_\la
\cap K_z) := \Pl; \ \ \ \T(\P_\la \cap B_{r_z}(z_0)) := \Plz.$$  We
also have the a.e. equality $B^{\la, z}=(r_z^d\la)^\b \B_{d-1}(\pi)
\times[0,(r_z^d\la)^{2\b})$, where $\B_{d-1}(\pi)$ is the closure of the
injectivity
region of  $\exp_{d-1}$. 

We next use the scaling transformations $\T$  on ${\mathcal A}_z(K)$
to define re-scaled $k$-face functionals $\xi^{\la, z}$ on re-scaled
point sets $\T(\P_\la \cap K_z)$; in the sequel we show that these
re-scaled functionals converge to the scaling limit functional
$\xi^{(\infty)}$ given in Definition \ref{xiinf}. In
the special case that $K$ is a ball, we remark that ${\mathcal
A}_z(K)= K$ for all $z \in \partial K$ and that $\T$ coincide for
all $z \in \partial K$, putting us in the set-up of \cite{CSY}.

\vskip.3cm

\noindent{\bf 3.3. Re-scaled $k$-face functionals $\txi, z \in
\partial K, \la \geq 1$}.
Fix $\la \in [1, \infty)$ and $z \in \pK$. 
{\em Let $\xi:= \xi_k$ be a generic $k$-face functional,} as given in Definition \ref{kface}.
 The inverse transformation
$[\T]^{-1}$ defines re-scaled $k$-face functionals $\txi(\w, \X)$
defined for $\w \in \K$ and $\X \subset \R^d$ by \be \label{txidef}
\txi(\w, \X):= \xi([\T]^{-1}(\w), [\T]^{-1}(\X \cap \K)). \ee It
follows for all $z \in \pK, \la \in [1, \infty),$ and $x \in K_z$
that
$
 \xi(x,
\P_\la \cap K_z):= \txi(\T(x), \Pl ).
$

We shall establish properties of the re-scaled $k$-face
functionals in the next section.  For now, we record
the distributional limit of the re-scaled  point processes $\Plz$ as $\la \to \infty$. 

\begin{lemm} \label{L3.3} Fix $z \in \pK$.  As $\la \to \infty$, we have $\Plz \tod \P$
in the sense of total variation convergence on compact sets.
\end{lemm}

{\em Proof.} This proof is a consequence of the discussion around
(2.14) of \cite{CSY}, but for the sake of completeness we include
the details. We  find the image by $\T$ of the measure on
$B_{r_z}(z_0)$ given by $\la r^{d-1} dr d\sigma_{d-1}(u)$. Under
$\T$ we have $h':= (r_z^d\la)^{2\b}(1 - {r \over r_z})$, whence $r = r_z(1
- (r_z^d\la)^{-2\b}h').$  Likewise we have $v':=(r_z^d\la)^{\b} v$, whence $v =
(r_z^d\la)^{-\b} v'$.
Under $\T$, the measure  $r^{d-1} dr$ becomes 
$$r^{d-1} dr = (r_z(1 - (r_z^d\la)^{-2 \b}))^{d-1}r_z^{1 - 2 \b d} \la^{-2 \b} dh'$$
and $d\sigma_{d-1}(u)$ transforms to
$$
d\sigma_{d-1}(u) = \frac{ \sin^{d-2}((r_z^d\la)^{-\b} |v'|) }{ |(r_z^d\la)^{-\b}
v' |^{d-2} } (r_z^d\la)^{- 1 + 2 \b}  dv'$$ as in (2.17) of \cite{CSY}.
Therefore the product measure $\la r^{d-1} dr d\sigma_{d-1}(u)$
transforms to \be \label{intensity} (1 - (r_z^d\la)^{-2 \b}h')^{d-1}
\frac{ \sin^{d-2}(\la^{-\b} |v'|) }{ |\la^{-\b} v' |^{d-2} } dh'
dv'.\ee
 The total variation distance between Poisson measures is upper
bounded by a multiple of the $L^1$ distance between their densities
(Theorem 3.2.2 in \cite{REI}) and since $(1 - (r_z^d\la)^{-2\b})^{(d-1)} \to 1$ as $\la \to \infty$, the result
follows.   \qed

\section{Properties of the re-scaled $k$-face functional  $\txi$
}

\allco

\noindent{\bf 4.1. Localization of  $\txi.$
}
We appeal to results of Reitzner \cite{Re} to show that
the re-scaled functionals $\txi$ given at \eqref{txidef} `localize', that is they are with
high probability determined by `nearby' point configurations.

For all $s > 0$  consider the inner parallel set of $\pK$, namely
\be \label{para} K(s):= \{x \in K \ : \ \de^{H}(x, \pK) \leq s
\},\ee with $\de^H$ being the Hausdorff distance.  Put \be \label{ee}  \e: = ( \frac{12 d
\log \la} { d_3 \la} )^{\b},\ee
where $d_3$ is as in Lemma 5 of Reitzner \cite{Re}.
Let $B_r(x)$ denote the Euclidean ball of radius $r$ centered at
$x$.  We begin with two localization properties  of the score $\xi$.
Here and elsewhere we shorthand $\xi_k$ by $\xi$.

\begin{lemm} \label{L4.0} (a) With probability at least $1 -
O(\la^{-4d})$, for all $z \in \pK$,  $\rho \geq 1$, we
have
 \be
\xi(x, \P_\la \cap K_z) = \begin{cases} \xi(x, \P_\la \cap
K_z(\rho \2e))
 & {\rm if} \ x \in
K_z(\2e)
\\
0 & {\rm if} \ x \in K_z \setminus K_z( \2e).
\end{cases}
 \ee

\noindent(b) There is a constant $D_1$ such that
 for all  $z \in \pK$, and $x \in K_z(\2e)$ we have
 $$
 P[\xi(x, \P_\la \cap K_z ) \neq \xi(x, \P_\la \cap K_z
 \cap B_{D_1 \e }(x))] = O(\la^{-4d}).
 $$

\end{lemm}

\noindent{\em Proof.} We prove part (a) with $\rho = 1$. The proof for $\rho > 1$ is identical.
 Let $X_i, i \geq 1,$ be i.i.d.
uniform on $K_z$. For every integer $l$, let $A_l$ be the event that
the boundary of $\rm{co}(X_1,...,X_l)$ is contained in
$K_z(\epsilon_l^2)$.

Following nearly verbatim the discussion on page 492 of  \cite{Re},
we note that $P[A_l^c]$ equals the probability that at least one
facet of $\rm{co}(X_1,...,X_l)$ contains a point distant at least
$\epsilon_l^2$ from the boundary of $K_z$,
i.e., this is the
probability that the hyperplane which is the affine hull of this
facet cuts off from $K_z$ a cap of height $\epsilon_l^2$ which
contains no point from $X_1,...,X_l$. By Lemma 5 of \cite{Re}, the
volume of this cap is bounded by $d_3 e_l^{d + 1} = 12d \log l/l$.

Thus  when $l$ is large enough so that   $(l-d)/l > 1/2$ (ie.  $l
> 2d$) and ${(12 d \log l)/ l} < 1$, and using $\log(1 -x) <
-x, \ 0 < x < 1$, we get \be \label{Re1} P[A_l^c] \leq \binom{l}{d}
\left( 1 - {12 d \log l \over l} \right)^{l - d}
 < l^d {1 \over d!} \exp\left( (l-d) (   -    {12 d \log l \over l}) \right)
  \leq  {l^d \over d!} l^{-6d}= {l^{-5d}\over d!}. \ee

Let $A_\la$ be the event that the boundary of $\rm{co}(\P_\la \cap
K)$ is contained in $K_z(\epsilon_\la^2)$. Letting $N(\la)$ be a
Poisson random variable with parameter $\la$ we compute
$$
P[A_\la^c] = \sum_{l= 0}^{\infty} P[A_l^c, N(\la) = l] < \sum_{| l -
\la | \leq \la^{3/4}} P[A_l^c] + P[|N(\la) - \la| \geq \la^{3/4}
].$$ The last term decays exponentially with $\la$ and so exhibits
growth $O(\la^{-4d})$.   By \eqref{Re1},  the first term has the
same growth bounds since
$$
\sum_{| l - \la | \leq \la^{3/4}} P[A_l^c] \leq 2 \la^{3/4} \max_{
|\la - l| \leq \la^{3/4} }P[A_l^c] \leq  2 \la^{3/4} {1 \over d!}
(\la - \la^{3/4})^{-5d} = O(\la^{-4d}),$$
concluding the proof of (a).


We prove assertion (b). By part (a), it suffices to show there is
$\rho_0 \geq 1$ such that for $x \in K_z(\2e)$
$$
 P[\xi(x, \P_\la \cap K_z(\rho_0 \2e) ) \neq \xi(x, \P_\la \cap K_z(\rho_0 \2e)
 \cap B_{D_1 \e }(x))] = O(\la^{-4d}).
 $$

We consider  the localization results described on pages 499-502 of
\cite{Re} and in the Appendix of \cite{Re}.  Using the set-up of
Lemma 6 of \cite{Re}, we choose  $m := m(\la):= \lfloor (d_6 \la/(4d + 1)
\log \la)^{(d-1)\b} \rfloor$ points $y_1,...,y_m$ on $\pK_z$ (here
$d_6$ is the constant of \cite{Re}) such that the Voronoi cells
$C_{\mbox{{\tiny Vor}}}(y_j), 1 \leq j \leq m,$ partition $K_z$, and such that the
diameter of $C_{\mbox{{\tiny Vor}}}(y_j) \cap \pK_z$ is $O( \e)$. Moreover, because all
$y_j$ are  on $\pK_z$, any bisecting line between two $y_j$ makes an
angle with $\partial K_z$ which is bounded from below.
Consequently, since the `width'
of $K_z( \2e)$ is $O(\2e)$, it follows that the diameter of the
truncated cells $C_{\mbox{{\tiny Vor}}}(y_j) \cap K_z(\2e)$ is also
$O( \e)$. Choose $\rho_0$  large enough so that $K_z(\rho_0 \2e)$
contains the caps $C_i, 1 \leq i \leq m,$ given near the end of page
498 of \cite{Re}.


For all $1 \leq j \leq m$, let $$S_j := \{ k \in \{1,2,...,m\}: \
C_{\mbox{{\tiny Vor}}}(y_k) \cap C(y_j, d_{10} m^{-2\b} ) \neq \emptyset \}$$ where $C(y,
h)$ denotes a cap at $y$ of height $h$, and where $d_{10}$ denote
the constant in \cite{Re}. Pages  498-500 of \cite{Re}
show  the
existence of a set $A^m$ such that $P[A^m] \geq 1 - c_{16}
\la^{-4d}$, and on $A^m$ the score $\xi(x, \P_\la \cap K_z(\rho_0
\2e))$ at  $x \in  K_z(\2e) \cap C_{\mbox{{\tiny Vor}}}(y_j)$ is determined by the
Poisson points belonging to \be \label{Re2} U_j := U_j(x):=  \bigcup_{k \in
S_j}
C_{\mbox{{\tiny Vor}}}(y_k) \cap K_z( \2e), \ee
where $j:= j(x) \in \{1,...,m\}$ is such that $C_{\mbox{{\tiny Vor}}}(y_j)$ contains $x$.
(Actually \cite{Re} shows this for the score $\xi(x, \P_\la\cap
K_z)$ and not for $\xi(x, \P_\la \cap K_z(\rho_0 \2e) )$, but the
proof is the same, since $\rho_0$ is chosen so that $K_z(\rho_0
\2e)$ contains the caps $C_i, 1 \leq i \leq m.$)
 By Lemma 7 of \cite{Re},  the cardinality of $S_j$
 is at most $d_8( d_{10}^{1/2} m^{-\b} m^\b  + 1)^{d + 1} =
O(1)$, uniformly in $1 \leq j \leq m$. 
 This implies
that on $A^m$, the score $\xi(x, \P_\la \cap K_z(\rho \2e))$ at  $x
\in  K_z(\2e) \cap C_{\mbox{{\tiny Vor}}}(y_j)$ is determined by the Poisson points in
$U_j$,   whose diameter is bounded by a constant multiple of the
diameter of the truncated cells $C_{\mbox{{\tiny Vor}}}(y_k)  \cap K_z(\rho \2e), k \in
S_j$, and is thus determined by points distant at most $D_1 \e$ from
$x$, $D_1$ a constant. Since $P[A_m^c] \leq c_{16} \la^{-4d}$, this
proves assertion (b). \qed


\ \

The next lemma shows
localization properties of $\txi$.  We first require more
terminology.


\begin{defn} For
all $z \in \pK$,
we put
 $$
S^{\la,z}:= \T( K_z(\2e)
\cap B_{2D_1 \e}(z)).$$
\end{defn}

Note that if $\w = (v',h') \in S^{\la,z}$, then $|v'| \leq D_2
(\log \la)^\b$  for some $D_2$ not depending on $z$ (here we use
$\sup_{z \in \partial K} r_z \leq C$).  Also, define $D_3$ by the relation
$2 [\sup_{z \in \partial K} r_z^{d \b} ] D_1
\la^\b \e = D_3(\log \la)^\b$. For all $L
> 0$ and $v \in \R^{d-1}$, denote by
 $\C_L(v)$  the cylinder $ \{(v',h) \in \R^{d-1} \times \R:
 \ |v' - v| \leq L \}$. Due to the non-linearity of  $T^{\la, z}$,
localization
 properties for $\xi$ do not in general imply localization properties
for  $\txi(\w, \Pl)$.  However, the next lemma 
says that if the
inverse image of $\w$ is close to $z$, then $\txi(\w, \Pl)$ suitably localizes.

\begin{lemm} \label{L4.1} Uniformly in
$z \in \pK$ and
$\w:=(v',h') \in S^{\la,z}$
we have
$$
P[ \txi(\w, \Pl ) \neq \txi(\w, \Pl   \cap \C_{D_3 (\log \la)^\b}
(v' )) ] = O (\la^{-4d}).
$$
\end{lemm}

\ \

\noindent{\em Remark.}  When $K$ is the unit ball we show in \cite{CSY} that the scores
$\txi$ localize in the following stronger sense: for all $\w :=
(v',h') \in K^{\la, z}$, there is an a.s. finite random variable
$R:= R(\w, \Pl)$ such that \be \label{local} \txi(\w, \Pl ) =
\txi(\w, \Pl \cap \C_{r} (v' ) )\ee for all $r \geq R$, with
$\sup_{\la} P[R>t]\to 0$ as $t \to \infty$.   We are
unable to show this latter property for arbitrary smooth $K$.

\ \

\noindent{\em Proof.} Fix the reference boundary point $z \in \pK$. Let $\rho_0$ be as in the proof of Lemma \ref{L4.0}(b).
For any $A \subset \R^+ \times \R^{d-1}$, we let $\T(A):= A^{\la, z}.$  In view
of Lemma \ref{L4.0}(b), it suffices to show for $\w:=(v',h') \in
S^{\la,z}$ that
$$
P[ \txi(\w, (\P_\la \cap K_z(\rho_0 \2e))^{\la, z} ) \neq \txi(\w,
(\P_\la \cap K_z(\rho_0 \2e))^{\la, z}  \cap \C_{D_3 (\log \la)^\b}
(v' )) ] = O(\la^{-4d}).
$$

Given $w'$, find $j:=j(w')$ such that $ C_{\mbox{{\tiny Vor}}}(y_j)$ contains $[\T]^{-1}(w'):=x$.
Recall the definition of $U_j:= U_j(x)$   at \eqref{Re2}
and recall that the proof of Lemma \ref{L4.0} shows that $\text{diam}(U_j) \leq D_1 \epsilon_\la
$. By the $C^3$ assumption, if $\la$ is large then  
 for all $z \in \pK$ the projection of
$U_j$ onto the osculating sphere at $z$ has a diameter comparable to
that of $U_j$, i.e., is generously bounded by $2D_1 \epsilon(\la)$.
Thus 
the spatial
diameter of $T^{\la,z}(U_j)$ is bounded by $  2 [\sup_{z \in \partial K} r_z^{d \b} ] \la^{\b} D_1 \e =
D_3(\log \la)^{\b}$, by definition of $D_3$. In other words
\be \label{4-4}
T^{\la,z}(U_j) \subset \C_{D_3 (\log \la)^\b}
(v' ).\ee


However, as seen in the proof of Lemma \ref{L4.0},
 with probability at least $1 - c_{16} \la^{-4d}$, the score
$\txi(\w, (\P_\la \cap K_z(\rho_0 \2e))^{\la, z} )$ is determined by
the points $(\P_\la \cap K_z(\rho_0 \2e))^{\la, z}$ in
$T^{\la,z}(U_j)$.
In view of \eqref{4-4}, the proof is complete.   \qed

\vskip.5cm

\noindent{\bf 4.2. Moment bounds for $\txi$.}
We use the localization results
to derive moment bounds for the re-scaled $k$-face functionals  $\txi$. For
a random variable $W$ and all $p > 0$, we let $||W||_p:= (\E
|W|^p)^{1/p}$.

\begin{lemm} \label{mom}
Let $\xi:= \xi_k, k \in \{0,1,...,d-1\}$.
 For all $p \in [1, 4]$ there are constants $M(p):= M(p,k) \in (0, \infty)$ such that
\be \label{first}\sup_{z \in \pK} \sup_{\la \geq 1} \sup_{\w \in
 B^{\la, z} } || \xi^{\la,z}(\w, \Plz) ||_p \leq M(p)
\ee and \be \label{second} \sup_{z \in \pK} \sup_{\la \geq 1}
\sup_{\w \in S^{\la,z}  } || \xi^{\la,z}(\w, \Pl)
||_p \leq M(p) (\log \la)^{k}. \ee
\end{lemm}

\noindent{\em Proof.}
The bound \eqref{first} follows as in Lemma 7.1 of
\cite{CSY}. To prove \eqref{second}, we argue as follows.
Given $z \in \pK$ and $\w \in S^{\la,z} $, we let
$$
E := E_z (\w):= \{ \txi(\w, \Pl ) = \txi (\w, \Pl   \cap \C_{D_3
(\log \la)^\b } (v' )  \cap (K_z(\2e))^{\la, z}    ) \}. $$  By
Lemmas \ref{L4.0}(a) and \ref{L4.1} we have $P[E^c] = O(\la^{-4d}).$

Let $N(s)$ be  a Poisson random variable with
parameter $s$.  
The cardinality of the point set
$$
\Pl   \cap \C_{D_3 (\log \la)^\b } (v' )  \cap (K_z(\2e))^{\la, z},
$$
is stochastically bounded by $N(C (\log \la)^{\b(d-1)} \cdot (\log
\la)^{2 \b}) = N(C \log \la)$, where $C$ is a generic constant whose
value may change from line to line.
 On the event $E$ the number of $k$-faces
containing $\w$ is generously bounded by $ \binom{ N(C \log \la) }
{k} \leq (N(C \log \la))^{k} $.

We now compute for $p \in [1,4]$:
$$
|| \xi^{\la, z}(\w, \Pl)||_p \leq || \xi^{\la, z}(\w,
\Pl){\bf{1}}(E)||_p + || \xi^{\la, z}(\w, \Pl){\bf{1}}(E^c)||_p. $$
The first term is bounded by $(k+1)^{-1}||N^{k}(C \log \la)||_p \leq
M(p) (\log \la)^{k}$.
The second term is bounded by
$$\frac{1}{k+1}
\left| \left| \binom {\text{card}(\P_\la)} {k} \right|
\right|_{pr} \la^{-4d/pq}, \ \ 1/r + 1/q = 1.
$$
We have $||  \binom {\text{card}(\P_\la)} {k} ||_{pr} =
O(\la^{k})$ and for $p \in [1,4]$  we may choose $q$ sufficiently
close to $1$ such that $\la^{-4d/pq} = O(\la^{-k})$. This gives
\eqref{second}.


\qed


\ \

{\em Remarks.} (i) Straightforward modifications of the proof of
Lemma \ref{L4.0} show that the $O(\la^{-4d})$ bounds of that lemma
may be replaced by $O(\la^{-md})$ bounds, $m$ an arbitrary integer,
provided that $\epsilon_\la$ given at \eqref{ee} is increased by a scalar multiple of
$m$.  In this way one could show that Lemma \ref{mom} holds for moments of any order $p > 0$.
 Since we do not require more than fourth moments for $\txi$, we
do not strive for this generality.

(ii) We do not claim that the bounds of Lemma \ref{L4.0}  are
optimal.  By McClullen's bound \cite{Mc}, the $k$ face
functional on an $n$ point set is bounded by $Cn^{d/2}$ and using
this bound for $k > d/2$ shows that the $(\log \la)^k$ term in
\eqref{second} can be improved to $(\log \la)^{d/2}$. The $\log \la$ factors could possibly be dispensed with
altogether, as mentioned in Section 5.3.

\vskip.5cm

\noindent{\bf 4.3. Comparison of scores for points in a ball and on
$K_z$.} The  $k$-face functional of Definition \ref{kface} on Poisson input on the ball is
well understood \cite{CSY}.  To exploit this we need to
show that the re-scaled functional $\txi$ on $\Pl$ is well
approximated by its value on $\Plz$.  We shall also need to show that
the pair correlation function for $\txi$ on $\Pl$ is well
approximated by the pair correlation function for $\txi$ on $\Plz$.
These approximations are established in the next four lemmas.

Our first lemma records a simple geometric fact. Locally
around $z$, the osculating ball to $K_z$ may lie inside or outside
$K_z$, but it is not far from $\pK_z$. The next lemma shows that the
distance decays like the cube of $|v'|$.

\begin{lemm} \label{L3.2} For all $z \in \pK$ and  $v:= (r_z^d \la)^{-\b} v'$ we have
\be \label{D4.4} r_z^{2\b d} \la^{2 \b}\left|1 - {f(\exp_{d-1} (v) )
\over r_z}\right| \leq D_4 r_z^{-1 - \b d}
\la^{-\b} |v'|^3. \ee 
\end{lemm}
{\em Proof.} We first show \eqref{D4.4}  when $d=2$. The boundary of
the osculating circle at $z$ coincides with $\pK$ up to at least
second order, giving $f(0) = r_z, f'(0) = f''(0) = 0$. The Taylor
expansion for $f$ around $0$ gives $|1 - {f(v) \over r_z}| \leq {1
\over 6} ||f^{'''}||_{\infty} r_z^{-1} |v|^3$, whence the result.

We now consider the case $d \geq 3$. Let $\exp_{d-1} (v) := \cos(|v|)
k_z + \sin(|v|) w$,
where $w := v/|v|$. 
It is enough to consider the section of the osculating ball
 and $K_z$ with the plane generated by $k_z$ and $w$. Indeed, we obtain in that
  plane a two-dimensional mother body with an osculating radius equal to $r_z$ at the point $z$.
  We may  apply the case $d=2$ to deduce the required result.  \qed

\ \

\begin{lemm} \label{L4.2} Uniformly for
  $z \in \pK$ and $\w \in S^{\la,z} \cap B^{\la, z}$, we have
\be \label{D4.2} \E \left| \txi(\w, \Pl) -  \txi(\w, \Plz) \right| =
O\left( 
 \la^{-\b/2} (\log \la)^{k + (\b + 1)/2}
\right). \ee
\end{lemm}
{\em Proof.}
For
$\w \in S^{\la,z} \cap B^{\la, z}$, we put \be \label{defE} E:=
E(\w) := \{ \txi(\w, \Pl ) = \txi(\w, \Pl \cap \C_{D_3 (\log
\la)^\b} (\w)) \} \ee
$$ \cup  \
\{ \txi(\w, \Plz) = \txi(\w, \Plz \cap \C_{D_3 (\log \la)^\b} (\w))
\},$$
so that $P[E^c] =
O(\la^{-4d})$ by Lemma \ref{L4.1}.
 Put
$$
F^{\la,z}(\w):= \txi(\w, \Pl) - \txi(\w, \Plz).$$ By Lemma \ref{mom}
with $p = 2$, we have $||F^{\la,z}(\w)||_2 \leq 2 M(2) (\log
\la)^{k},$ uniformly in $\w, \la$ and $z$.

 Recall $\w := (v',h')$.
 For all $\w \in S^{\la,z} \cap B^{\la, z}$
put \be \label{reg} R(\w):= \{(v'',h'') \in \R^{d-1} \times \R : \ |v'' -
v'| \leq D_3 (\log \la)^\b, $$ $$ \ \ \ \ \ |h''| \leq (r_z^d \la)^{2\b}|1 - r_z^{-1}f(
\exp_{d-1}((r_z^d \la)^{-\b} v'') ) |\}. \ee Write
$$
\E |F^{\la,z}(\w)| = \E |(F^{\la,z}(\w))({\bf 1}(E) + {\bf 1}(E^c))|
.$$ On $E$ we have $F^{\la,z}(\w)= 0$, unless the realization of
$\Pl$ puts points in the set $R(\w)$.  By the Cauchy-Schwarz inequality and Lemma \ref{mom} with $p = 2$ there,  we have \be
\label{CS} \E |( F^{\la,z}(\w)){\bf 1}(E)| \leq 2M(2) (\log \la)^k
(P[ {\bf 1}(\Pl \cap R(v') \neq \emptyset)])^{1/2}. \ee The Lebesgue
measure of $R(\w)$ is bounded by the product of  the area of its
`base', that is $(2D_3 (\log \la)^\b)^{d-1}$ and its `height', which by Lemma
\ref{L3.2} is at most $D_4 r_z^{-1- \b d}\la^{-\b} \left(|v'| + D_3 (\log
\la)^\b\right)^3.$
By \eqref{intensity}, the $\Pl$ intensity measure
of $R(\w)$, denoted by $|R(\w)|$,  thus satisfies
\be \label{Rbd} |R(\w)| \leq (2D_3 (\log \la)^\b)^{d-1}D_4 r_z^{-1- \b d}\la^{-\b}
(|v'| + D_3 (\log \la)^\b)^3. \ee
 Since $1 -
e^{-x} \leq x$ holds for all $x$
it follows that \be \label{L4.3a}
P[ {\bf 1}(\Pl \cap R(v') \neq \emptyset)] = 1 - \exp(-|R(\w)| ) \leq |R(\w)|.
\ee
Combining \eqref{CS}-\eqref{L4.3a}, and recalling that $|v'| \leq D_2 (\log \la)^\b$, shows that $\E
|(F^{\la,z}(\w)){\bf 1}(E)|$ is bounded by the right hand side of
\eqref{D4.2}.

 Similarly, Lemma \ref{mom}, the bound $P[E^c] = O( \la^{-4d})$, and the Cauchy-Schwarz inequality  give
  $\E|(F^{\la,z}(\w)){\bf
1}(E^c)|
= O((\log \la)^k \la^{-2d})$, which is dominated by the right
hand side of \eqref{D4.2}. Thus
\eqref{D4.2} holds as claimed.
 \qed

\ \

The next lemma is the analog of Lemma 7.2 in \cite{CSY}. It justifies the
use of the scaling limit terminology for $\xi^{(\infty)}$, as given by
Definition \ref{xiinf}.

\begin{lemm} \label{L4.3} For all $z \in \pK$ and $(\0,h) \in
\K$ we have
$$
\lim_{\la \to \infty}| \E \txi((\0,h), \Pl) - \E
\xi^{(\infty)}((\0,h), \P) |
= 0.
$$
\end{lemm}

\noindent{\em Proof}. We bound $| \E \txi((\0,h),
\Pl) - \E \xi^{(\infty)}((\0,h), \P) |$ by
$$ | \E \txi((\0,h), \Pl) - \E \txi((\0,h), \Plz) | + | \E
\txi((\0,h), \Plz) - \E \xi^{(\infty)}((\0,h), \P)|.$$ The
first term goes to zero by Lemma \ref{L4.2} with $w' = (\0,h)$ and
the second term goes to zero by Lemma 7.2 of \cite{CSY}. \qed

\ \

We next recall the definition of the pair correlation function for
the score $\xi$ as well as for its re-scaled version.

\begin{defn} (Pair correlation functions) For all $x, y \in K_z$,
any random point set $\Xi \subset K_z$, and any $\xi$  we put \be
\label{pcorr1} c(x,y; \Xi) := c^\xi(x,y; \Xi):= \E \xi(x, \Xi \cup
y) \xi(y, \Xi \cup x) - \E \xi(x, \Xi ) \E \xi(y, \Xi). \ee For all
$\la \geq 1, z \in \pK, \ (\0, h) \in \K,$ and $(v',h') \in \K$,
define  the re-scaled pair correlation function of the $k$-face
functional as \be \label{pcorr2} \c((\0,h), (v',h'); \Pl):=
$$ $$ \E \txi((\0,h), \Pl \cup (v',h')) \txi((v',h'), \Pl \cup (\0,h)) -
\E \txi((\0,h), \Pl ) \E \txi((v',h'), \Pl). \ee
\end{defn}

The next lemma shows that the pair correlation function for $\txi$
on $\Pl$ is well approximated by the pair correlation function for
$\txi$ on $\Plz$.

\begin{lemm} \label{L4.5}  Uniformly for
 $z \in \pK$, $\wo := (\0, h) \in S^{\la,z}
\cap B^{\la, z}$ and $ \w:= (v',h') \in S^{\la,z} \cap B^{\la,
z}$, we have \be \label{D4.5} |\c(\wo, \w; \Pl) - \c(\wo, \w; \Plz)|
= O \left( 
\la^{-\b/3}(\log \la)^{2k + \b(d + 2)/3}
\right). \ee
\end{lemm}

\noindent{\em Proof.}  It suffices to modify the proof of Lemma \ref{L4.2}.
Put  $F:= E(\wo) \cap E(\w)$, where $E(\wo)$ and $E(\w)$ are defined
at \eqref{defE}. We have $P[F^c] = O( \la^{-4d})$ by Lemma
\ref{L4.1}. Write
$$
\E \txi(\wo, \Pl \cup (v',h')) \txi(\w, \Pl \cup \wo)  - \E
\txi(\wo, \Plz \cup \w) \txi(\w, \Plz \cup \wo)$$ $$
  = \E[ \{\txi(\wo, \Pl \cup \w) \txi(\w, \Pl \cup
\wo)   - \txi(\wo, \Plz \cup \w) \txi(\w, \Plz \cup \wo)\} {\bf{1}}
(F)] \ + \ $$
 \be \label{D1}
+ \ \E[\{\txi(\wo, \Pl \cup \w) \txi(\w, \Pl \cup \wo)  - \txi(\wo,
\Plz \cup \w) \txi(\w, \Plz \cup \wo)\} {\bf{1}}(F^c)]
\ee
$$
:= I_1 + I_2.
$$

The random variable  in the  expectation $I_1$ vanishes, except on
the event
$$H(\wo,\w):= \{ \Pl \cap R(\w) \neq \emptyset \} \cup \{ \Pl \cap R(\wo) \neq \emptyset \},$$
where $R(\w)$ and $R(\wo)$ are at \eqref{reg}.
The H\"older inequality $||UVW||_1 \leq ||U||_{3} ||V||_{3} ||W||_3$
for random variables $U, V, W$ and Lemma \ref{mom} imply that
$$I_1 \leq 2 (M(3))^2 (\log \la)^{2k} (P[H(\wo,\w)])^{1/3},$$ that is to say
$$
I_1 = O \left(  (\log \la)^{2k} \left(  r_z^{-1- \b d}\la^{-\b}
(\log \la)^{\b(d-1)} [ (|v'| + D_3 (\log \la)^\b )^3 + (D_3 (\log
\la)^\b)^3 ] \right)^{1/3} \right),
$$
which for $|v'| \leq D_2(\log \la)^\b$
satisfies the growth bounds on  the right hand
side of \eqref{D4.5}. 

Now term $I_2$ in   \eqref{D1} is bounded by $2(M(3))^2 (\log
\la)^{2k}(P[F^c])^{1/3}$, which is of smaller order than the
right hand side of \eqref{D4.5}.  This shows that \eqref{D1} also
satisfies the growth bounds on  the right hand side of \eqref{D4.5}.


It remains to bound \be \label{D2} |\E \txi(\wo, \Pl ) \E \txi(\w,
\Pl ) - \E \txi(\wo, \Plz ) \E \txi(\w, \Plz )|.\ee   Notice that
the difference \eqref{D2} differs from \be \label{D3} |\E \txi(\wo,
\Pl ){\bf 1}(F) \E \txi(\w, \Pl ){\bf 1}(F)  - \E \txi(\wo, \Plz
){\bf 1}(F) \E \txi(\w, \Plz ){\bf 1}(F)|\ee by at most \be
\label{D3a}4(M(3))^2 (\log \la)^{2k} (P[F^c])^{1/3} \leq
C(M(3))^2 (\log \la)^{2k}  \la^{-4d/3}, \ee which is of
smaller order than  the right hand side of \eqref{D4.5}.

Now we control the difference \eqref{D3} which we write as $|\E e_1
\E e_2 - \E e_3 \E e_4|$, where $e_1 := \txi(\wo, \Pl ){\bf 1}(F)$,
$e_2:= \txi(\w, \Pl ){\bf 1}(F)$, $e_3:= \txi(\wo, \Plz ){\bf
1}(F)$, and $e_4:= \txi(\w, \Plz ){\bf 1}(F)$. The proof of Lemma
\ref{L4.2} (with $E$ replaced by $F$) shows that \be \label{D4} \E
|e_1 - e_3| = O( 
\la^{-\b/2 } (\log \la)^{k + (\b +
1)/2}) \ee and \be \label{D5} \E |e_2 - e_4| =  O( 
\la^{-\b/2 } (\log \la)^{k + (\b + 1)/2}) \ee Since $ |\E e_1 \E
e_2 - \E e_3 \E e_4| \leq |\E e_1| |\E e_2 - \E e_4| + |\E e_4 | |\E
e_1 - \E e_3| $ it follows that \eqref{D2} is bounded by \be
\label{L4.6a} O( 
\la^{-\b/2 } (\log \la)^{2k +
(\b + 1)/2})  + O( (\log \la)^{2k}  \la^{-4d/3}), \ee i.e., is
bounded by the right-hand side of \eqref{D4.5}.
 \qed

\ \

Our last lemma describes a decay rate for $c(x,y; \P_\la \cap
K_z)$, a technical fact used in the sequel.

\begin{lemm} \label{L4.4} For all $z \in \pK$ and $x ,y \in K_z(\2e)$  with $|x - y| \geq 2D_1 \e$, we have
$$\lim_{\la \to \infty} \la^{ 1 + 2 \b } c(x,y; \P_\la \cap K_z) =0.
$$
\end{lemm}

\noindent{\em Proof.} Fix $x \in K_z(\2e)$. {\em To lighten the notation we
abbreviate $\P_\la \cap K_z$ by $\P_\la$ in this proof only.} For $y
\in K_z(\2e)$, put
$$
E:= E(x,y):= \{\xi(x, \P_\la) = \xi(x, \P_\la \cap B_{D_1 \e}(x)) \}
\cup \{\xi(y, \P_\la) = \xi(y, \P_\la \cap B_{D_1 \e}(y)) \}.
$$
Lemma \ref{L4.0}(b) gives
 \be \label{exp2}  P[E^c]
= O( \la^{-4d}). \ee If $|x - y| \geq 2D_1 \e$, then $\xi(x, \P_\la
\cup y)$ and $\xi(y, \P_\la \cup x)$ are independent on $E$, giving
$$\E [\xi(x, \P_\la \cup y) \xi(y, \P_\la \cup x) {\bf 1}(E)] = \E
[\xi(x, \P_\la) {\bf 1}(E) \cdot \xi(y, \P_\la) {\bf 1}(E)]$$
$$ = \E [\xi(x, \P_\la) {\bf 1}(E)] \cdot \E [ \xi(y,\P_\la)
{\bf 1}(E)].$$ Writing ${\bf 1}(E) = 1 - {\bf 1}(E^c)$ gives
$$\E [\xi(x, \P_\la \cup y) \xi(y, \P_\la \cup x)
{\bf 1}(E)]$$
$$ =  \left(\E \xi(x, \P_\la ) -  \E [\xi(x, \P_\la ) {\bf
1}(E^c)]\right) \cdot  \left(\E \xi(y, \P_\la) -  \E [\xi(y,\P_\la )
{\bf 1}(E^c)]\right) $$ $$ = \E \xi(x, \P_\la)\E \xi(y, \P_\la) +
G(x,y),$$ where
$$G(x,y) : = - \E \xi(x, \P_\la )\E [\xi(y, \P_\la ) {\bf
1}(E^c)]$$
$$
- \ \E\xi(y, \P_\la)\E [\xi(x, \P_\la ) {\bf 1}(E^c)] +  \E
[\xi(x,\P_\la ) {\bf 1}(E^c)] \cdot \E [\xi(y, \P_\la ) {\bf
1}(E^c)].$$

Let $N(\la) := \text{card}(\P_\la \cap K_z)$.  By McClullen's bounds
\cite{Mc} for the number of $k$-dimensional faces and standard
moment bounds for Poisson random variables we have $||\xi(x,
\P_\la)||_1 \leq C ||N^{d/2}(\la)||_1 \leq C \la^{d/2}$ and similarly
$ ||\xi(y, \P_\la)||_2 \leq C \la^{d/2}.$ By the Cauchy-Schwarz inequality, it
follows that
$$|\E \xi(x, \P_\la )\E [\xi(y, \P_\la ) {\bf 1}(E^c)]| =
O(\la^{d/2} \la^{d/2} (P[E^c])^{1/2}) = o(\la^{-1 - 2\b }),$$
where the last estimate easily follows by \eqref{exp2}.
 The
other two terms comprising $G(x,y)$ have the same asymptotic
behavior and so $G(x,y) = o(\la^{-1 - 2\b })$.


On the other hand, $\E [\xi(x, \P_\la \cup y) \xi(y, \P_\la \cup x)
{\bf 1}(E)]$ differs from $\E \xi(x, \P_\la \cup y)) \xi(y, \P_\la
\cup x)$ by $\E [\xi(x, \P_\la \cup y) \xi(y,\P_\la \cup x) {\bf
1}(E^c)].$  The H\"older inequality $||U V W||_1 \leq ||U||_4 ||V||_4 ||W||_2$
 shows  that this term is $o(\la^{-1 - 2\b })$.


Thus $\E \xi(x, \P_\la \cup y)) \xi(y,\P_\la \cup x)$ and $\E \xi(x,
\P_\la )\E \xi(y, \P_\la)$ differ from $\E [\xi(x, \P_\la \cup y)
\xi(y,\P_\la \cup x) {\bf 1}(E)]$ by $o(\la^{-1 - 2\b })$, concluding the proof of
Lemma \ref{L4.4}. \qed

\section{Proof of Theorem \ref{Th2}}

\allco

 Recall that  ${\mathcal M}(K)$ denotes the medial axis
 of $K$ and, for every $z\in \partial K$   the inner unit-normal vector
 of $\partial K$ at $z$ is $k_z$.  Put $t_{m}(z):=\inf\{t>0:z+ tk_z\in {\mathcal M}(K)\}$.
 Then  the map $\varphi:(z,t) \longmapsto (z+tk_z)$ is a
  diffeomorphism from $\{(z,t):z\in \partial K,0<t<t_m(z)\}$ to
  $\mbox{Int}(K)\setminus {\mathcal M}(K)$. In particular, $z\longmapsto -k_z$ is
  the Gauss map and its differential is the shape operator or Weingarten map
  $W_z$, which we recall has eigenvalues $C_{z,1},\cdots,C_{z,d-1}$.
  Consequently, the Jacobian of $\varphi$ may be written as $\det(I-tW_z)=\prod_{i=1}^{d-1}(1-tC_{z,i}).$ 

\ \

\noindent{\bf 5.1. Proof of expectation asymptotics \eqref{main3}.}
 Fix $g \in
\C(K)$ and let $\xi$ and $\mu_\la^\xi$ denote a generic $k$ face functional and $k$ face measure,
respectively. Recall that we may uniquely write $x \in K \setminus {\mathcal M}(K)$
as $x := (z,t)$, where $z \in \partial K$, and $t \in (0, \infty)$ is the
distance between $x$ and $z$.  Write
$$
\la^{-1 + 2\b} \E [\langle g, \mu_\la^\xi \rangle] = \la^{2\b } \int_K
g(x) \E \xi(x, \P_\la \cap K) dx
$$
$$
= \la^{2\b} \int_{z \in \pK} \int_0^{t_m(z)} g((z,t)) \E \xi((z,t),
\P_\la \cap K) \cdot \Pi_{i=1}^{d-1} (1 - t C_{z,i}) dt dz.
$$
For each $z \in \pK$, we apply the transformation ${\cal A}_z$ to
$K$.  Recalling from \eqref{stable} that $\xi$ is stable under
${\cal A}_z$, we have $\E \xi((z,t), \P_\la \cap K) = \E \xi((z,t),
\P_\la \cap K_z)$, since ${\cal A}_z(z, t):=(z, t)$  and
${\cal A}_z(\P_\la \cap K) \eqd \P_\la \cap K_z$.  It follows that
$$
\la^{-1 + 2\b} \E [\langle g, \mu_\la^\xi \rangle] = \la^{2\b}
\int_{z \in \pK} \int_0^{t_m(z)} g((z,t)) \E \xi((z,t), \P_\la \cap
K_z)  \cdot \Pi_{i=1}^{d-1} (1 - t C_{z,i}) dt dz.
$$
By Lemma \ref{L4.0}(a), the bound \eqref{second} with $p = 2$, and the
Cauchy-Schwarz inequality, it follows that uniformly in $x \in K_z
\setminus K_z(\2e)$ we have $\lim_{\la \to \infty} \la^{2 \b} \E
\xi(x, \P_\la \cap K_z) = 0$.
Since $$\sup_{\la \geq 1} \sup_{x \in
K_z \setminus K_z(\2e) } \la^{2 \b} \E \xi(x, \P_\la \cap K_z) \leq
C,$$ the bounded convergence theorem shows that we can restrict the
range of integration of $t$ to the interval  $[0, \2e]$ with error
$o(1)$. This gives \be \label{intdom} \la^{-1 + 2\b} \E [\langle g,
\mu_\la^\xi \rangle] = \la^{2\b} \int_{z \in \pK} \int_0^{\2e}
g((z,t)) \E \xi((z,t), \P_\la \cap K_z ) \cdot \Pi_{i=1}^{d-1} (1 -
t C_{z,i}) dt dz + o(1). \ee

Changing variables with $t = r_z (r_z^d \la)^{-2 \b}h$ and using $h
= (r_z^d \la)^{2 \b} (r_z - r)/r_z = (r_z^d \la)^{2 \b}(t/r_z)$
gives $\xi((z,t), \P_\la \cap K_z) =  \xi^{\la, z}((\0, h), \Pl )$.
Letting $h(\la, z):= r_z^{-1 + 2 \b d} \la^{2 \b} \epsilon_\la^2$ we
get
$$
\la^{-1 + 2\b} \E [\langle g, \mu_\la^\xi \rangle]
$$
$$
 = \int_{z \in \pK} r_z^{1 - 2 \b d} \int_0^{h(\la, z)  } g((z, o_u(1))) \E
\xi^{\la, z}((\0, h), \Pl) \cdot \Pi_{i=1}^{d-1} (1 - o_u(1) )  dh
dz + o(1)
$$
where  $o_u(1)$ denotes a quantity tending to zero as
$\la \to \infty$, uniformly in $z \in \partial K$ and uniformly in $h \in [0, h(\la, z)]$, not
 necessarily the same at each occurrence.%

Note that $(\0,h)$ belongs to $S^{\la, z} \cap B^{\la, z}$ and so we may apply Lemma \ref{L4.2}
to $\xi^{\la, z}((\0, h), \Pl)$.   Thus, with $\w$
set to $(\0,h)$ in Lemma \ref{L4.2}, we have
$$
\sup_{z \in \partial K} \sup_{h \in [0, h(\la, z)]} h(\la, z) \left| \E \txi((\0, h), \Pl) -
\E \txi((\0, h), \Plz) \right| = o(1),
$$
and so we may replace
$\E \txi((\0, h), \Pl)$ by $\E \txi((\0, h), \Plz)$ with error
$o(1)$. 
We also have
$r_z^{1 - 2 \b d} = \ka(z)^{1/(d + 1)}$.
In other words,
$$
\la^{-1 + 2\b} \E [\langle g, \mu_\la^\xi \rangle] $$ $$ = \int_{z
\in \pK} \ka(z)^{1/(d + 1)} \int_0^{h(\la, z) } g((z, o_u(1)) \E
\txi((\0, h), \Plz) \cdot \Pi_{i=1}^{d-1} (1 - o_u(1) ) dh dz +
o(1).
$$
By Lemma 3.2 of \cite{CSY}, the integrand is dominated by an
exponentially decaying function of $h$, uniformly in $z$ and $\la$.

The continuity of $g$, 
and the dominated convergence theorem give
\begin{equation}\label{limitexpinterm}
\lim_{\la \to \infty} \la^{-1 + 2\b} \E [\langle g, \mu_\la^\xi
\rangle] = \int_{z \in \pK}  g(z)  \ka(z)^{1/(d + 1)} \left[ \int_0^{\infty} \E
\xi^{(\infty)}((\0, h), \P)  dh  \right] dz.
\end{equation}
This gives \eqref{main3}, as desired.

\ \

\noindent{\bf 5.2. Proof of variance asymptotics \eqref{main4}. }
Recalling \eqref{pcorr1}, for fixed $g \in \C(K)$ we have
$$
\la^{-1 + 2\b} \Var[\langle g, \mu_\la^\xi \rangle] = \la^{2\b}
\int_K g(x)^2 \E \xi^2(x, \P_\la \cap K) dx \ +
$$
$$
+ \ \la^{1+ 2\b}  \int_{K} \int_{K} g(x) g(y) c(x,y; \P_\la \cap K)  dy dx :=
I_1(\la) + I_2(\la).
$$
Following the proof of  \eqref{main3} until  \eqref{limitexpinterm} shows that  \be \label{V1} \lim_{\la \to
\infty}I_1(\la) = \int_{z \in \pK}  g(z)^2  \ka(z)^{1/(d + 1)}
\int_0^{\infty} \E (\xi^{(\infty)}((\0, h), \P))^2  dh dz. \ee

Turning to  $I_2(\la)$, 
write $x$ in curvilinear coordinates $(z,t)$ with respect to $\pK$.
This gives $dx = \Pi_{i=1}^{d-1}(1 - t C_{z,i})dt dz$.
Apply the map ${\cal A}_z$, write
${\cal A}_z(y) = \y$ for $y \in K$,  and use stability
\eqref{stable} to get \be \label{double} I_2(\la) = \la^{1 + 2\b}
\int_{z \in \pK} \int_0^{t'_m(z)}  \int_{\y \in K_z }  g((z,t))
g(\y) c((z,t), \y; \P_\la \cap K_z) d\y \cdot \Pi_{i=1}^{d-1}(1 - t
C_{z,i}) dt dz + o(1). \ee Here
$$
c((z,t), \y; \P_\la \cap K_z) = \E\xi((z,t), \P_\la \cap K_z \cup \y)
\xi(\y, \P_\la \cap K_z \cup (z,t)) -
 \E\xi((z,t), \P_\la \cap K_z) \E \xi(\y, \P_\la \cap K_z ).
$$
The McClullen bound \cite{Mc} gives  \be \label{cbd}
|c((z,t), \y; \P_\la \cap K_z)| \leq C\E [N(\la)^d] \leq C \la^d, \ee
where here $N(\la)$ denotes the cardinality of $\P_\la \cap K_z$.

We make the following three modifications to the triple integral
\eqref{double}, each one giving an error of $o(1)$:

(i) Replace the integration domain $\{\y \in K_z \}$ by $\{\y \in
K_z(\2e)\}$.  
 Indeed,
uniformly in $\y \in K_z \setminus K_z(\2e)$ we have
$$
\lim_{\la \to \infty} \la^{1 + 2 \b} c((z,t), \y; \P_\la \cap K_z) =
0,$$ by Lemma \ref{L4.0}(a), the bound \eqref{cbd}, and
the Cauchy-Schwarz inequality.  Since
$$
\sup_{\la \geq 1} \sup_{(z,t), \y \in K_z \setminus K_z(\2e)}
 \la^{1 + 2 \b} c((z,t), \y; \P_\la \cap K_z) \leq C,$$ the assertion follow by  the bounded
convergence theorem.


(ii) Replace the integration domain $\{\y \in K_z(\2e)\}$ by $\{\y \in
K_z( \2e) \cap B_{2D_1 \e}((z,t))\}$ (use  Lemma \ref{L4.4} and the
bounded convergence theorem). 

(iii) Replace the integration domain $[0,t'_m(z)]$ by $[0,\2e]$, as
at \eqref{intdom}.

These modifications yield   \be \label{double1}
I_2(\la) = \ee
 $$= \la^{1 + 2\b} \int_{z \in \pK} \int_0^{\2e}  \int_{\y \in K_z(\2e) \cap B_{2D_1 \e}((z,t)) }
  g((z,t)) g(\y) c((z,t), \y; \P_\la \cap K_z) d\y \cdot \Pi_{i=1}^{d-1}(1 - t C_{z,i}) dt dz + o(1).
$$
\ \

Changing variables with $\y = (r,u)$ gives $d\y = r^{d-1}dr
d\sigma_{d-1}(u)$ and it also gives
$$
\T((r,u)) = ( (r_z^d \la)^\b \exp_{d-1}^{-1}(u), (r_z^d \la)^{2 \b}(1 - {r \over
r_z})) = ((r_z^d \la)^\b v, h') = (v',h') = \w.
$$
Thus the covariance $c((z,t), \y; \P_\la \cap K_z)$ transforms to
$\c((\0,h), (v',h'); \Pl)$. Now change variables with  $t = r_z
(r_z^d \la)^{-2\b}h, v' = (r_z^d \la)^{\b} v, $ and $h' = (r_z^d \la)^{2\b}(1 -
\frac{r}{r_z})$.

The differential $\la^{1 + 2 \b}\Pi_{i=1}^{d-1}(1 - t C_{z,i})
r^{d-1}dr d \sigma_{d-1}(u) dt dz$  transforms to the differential
$$
\la^{1 + 2 \b} \Pi_{i=1}^{d-1}(1 - r_z (r_z^d \la)^{-2 \b} h C_{z,i}) ((1 -
(r_z^d \la)^{-2\b }h')r_z)^{d-1} r_z (r_z^d \la)^{-2 \b} dh' \times (r_z^d \la)^{-\b (d-1)}
dv' r_z (r_z^d \la)^{-2\b} dh dz$$
$$
= \Pi_{i=1}^{d-1}(1 - r_z (r_z^d \la)^{-2 \b} h C_{z,i}) (1 -
(r_z^d \la)^{-2\b }h')^{d-1} r_z^{1 - 2 \b d} dh' dv' dh dz.
$$

The upper limit of integration $\2e$ in \eqref{double1} changes to
$h(\la, z) $ and the domain of integration $K_z(\2e) \cap B_{2D_1
\e}((z,t))$ gets mapped to $S^{\la, z}$.   This gives \be \label{V2}
I_2(\la) = \int_{z \in \pK} \int_0^{h(\la, z)} \int_{(v',h') \in
S^{\la,z} } G_{\la}(h',v',h, z) dh'dv' dh dz + o(1),\ee where,
recalling $r_z^{1 - 2\b d} = \ka(z)^{1/(d + 1)}$, we get
$$G_{\la}(h',v', h, z) :=
 \ka(z)^{1/(d + 1)}
g((z, o_u(1))) g(r_z(1- o_u(1)), (r_z^d \la)^{-\b}v')  $$ $$ \cdot \c((\0,h), (v',h'); \Pl)
 \Pi_{i=1}^{d-1}(1 - o_u(1)) (1 -
o_u(1))^{d-1}.$$
We next restrict the integration domain $S^{\la,z}$ to $S^{\la,z} \cap B^{\la, z}$ since by Lemma \ref{L3.2} and
the moment bounds \eqref{second} we have
$$
\int_{z \in \pK} \int_0^{h(\la, z)} \int_{(v',h') \in
S^{\la,z}  \cap (B^{\la, z})^c } G_{\la}(h',v',h, z) dh'dv' dh dz = o(1).
$$
By Lemma \ref{L4.5}, uniformly on the range $\{(v',h') \in
S^{\la,z}  \cap B^{\la, z} \}$ and uniformly over $h \in [0, h(\la,z)]$, the covariance term $\c((\0,h), (v',h'); \Pl)$
differs from the covariance term $\c((\0,h), (v',h'); \Plz)$ by a
term of order $ \la^{-\b/3}$, modulo logarithmic terms. The integral
of this difference over $$(h',v',h,z) \in S^{\la,z}  \times [0,
h(\la, z)] \times \pK$$ is also $o(1).$ This gives \be \label{V4}
I_2(\la) = \int_{z \in \pK} \int_{|h| \leq h(\la, z)} \int_{(v',h')
\in S^{\la,z}  \cap B^{\la, z} } \tilde{G}_{\la}(h',v',h, z) dh'dv' dh dz + o(1),\ee
where
$$
 \tilde{G}_{\la}(h', v', h, z)
= \ka(z)^{1/(d + 1)} g((z, o_u(1) ) g(r_z(1 - o_u(1), o_u(1) ) $$ $$ \cdot \c((\0,h), (v',h'); \Plz)
\Pi_{i=1}^{d-1} (1 - o_u(1)) (1 -
o_u(1))^{d-1}.$$

Recalling the definition of $\zeta_{\xi^{(\infty)}}$ at \eqref{SO2}
we get via Lemma 7.2 of \cite{CSY} that
$$ \lim_{\la \to \infty} \tilde{G}_{\la}(h', v', h, z)
 = \ka(z)^{1/(d + 1)} g(z)^2   \ \zeta_{\xi^{(\infty)}}((\0,h), (v',h'); \P).$$

The first part of Lemma 7.3 of \cite{CSY} shows  that
$\c((\0,h), (v',h'); \Plz)$ is dominated  by an integrable function of
$h',v',h$ and $z$ on $[0, \infty)  \times \R^{d-1} \times [0,
\infty) \times \pK$.
 Since $\sup_{z \in \pK} | \ r_z^{d+1}|$ and $||g||_{\infty}$ are both
bounded and since the integration domain $S^{\la,z}  \cap B^{\la, z} $ increases up
to $\R^{d-1} \times [0, \infty)$, the dominated convergence theorem
gives \be \label{V5} \lim_{\la \to \infty} I_2(\la) = \int_{z \in
\pK} g(z)^2 \ka(z)^{1/(d + 1)}  \int_0^{\infty} \int_{\R^{d-1}}
\int_0^{\infty} \zeta_{\xi^{(\infty)}}((\0,h), (v',h'); \P) dh'
dv' dh dz. \ee Combining \eqref{V1} and \eqref{V5} gives \be
\label{V6} \lim_{\la \to \infty} \la^{-1 + 2 \b} \Var[\langle g,
\mu_\la^\xi \rangle] = \int_{\pK} g(z)^2  \ka(z)^{1/(d + 1)}
\int_0^{\infty}  \E (\xi^{(\infty)}((\0, h), \P))^2  dh \ + $$ $$ +
\ \ka(z)^{1/(d + 1)} \int_0^{\infty} \int_{\R^{d-1}} \int_0^{\infty}
\zeta_{\xi^{(\infty)}}((\0,h), (v',h'); \P) dh' dv' dh  dz. \ee



Recalling the definition of  $\sigma^2(\xi^{(\infty)})$ at
\eqref{S03}, this yields \be \label{V9} \lim_{\la \to \infty}
\la^{-1 + 2 \b} \Var[\langle g, \mu_\la^\xi \rangle] =
\sigma^2(\xi^{(\infty)}) \int_{\pK} g(z)^2 \ka(z)^{1/(d + 1)} dz.
\ee This concludes the proof of variance asymptotics and the proof
of Theorem \ref{Th2}. \qed

\vskip.5cm

\noindent{\em Remark.}
If one could show that $\txi$ localize in the sense of \eqref{local},
then one could show that the moment bounds of Lemma \ref{mom} are
independent of $\la$.   We expect that one could subsequently  weaken
the $C^3$ boundary assumption to a $C^2$ assumption by making these three
changes: (i) replace the right-hand side of \eqref{D4.4} with
$o(1)|v'|^2$, (ii) in Lemmas \ref{L4.2} and \ref{L4.5}, drop the
restrictions $\w_0, \w \in S^{\la,z} \cap B^{\la, z},$  and  replace
the bounds on the right-hand side of \eqref{D4.2} and \eqref{D4.5}
with $o(1)$ bounds, and (iii) show that $ \c((\0,h), (v',h'); \Pl)$
decays exponentially in $|v'|$ and $h'$, showing that $G_\la(h',v',
h, z)$ is integrable. We could then  directly apply the dominated
convergence theorem to $\E \xi^{\la, z}((\0, h), \Pl)$ and
$\c((\0,h), (v',h'); \Pl)$ without needing the error approximations
of Lemmas \ref{L4.2} and \ref{L4.5}.
\section{Proof of Theorem \ref{Th1b}}
\allco
The image of $K$ by $x\longmapsto {\rm{vol}}(K)^{-1/d}\cdot x$ is a convex
 body of unit volume so without loss of generality, we may assume in this section that ${\rm{vol}}(K)=1$.
  The proof of Theorem \ref{Th1b} via Theorem \ref{Th1} is a rewriting of a result previously obtained
  by Vu (see \cite{Vu}, Proposition 8.1) in the case $k=0$. For sake of completeness, we include here a
  proof which does not use any large deviation result for $f_k(K_\la)$. The method uses  a coupling of
   the Poisson point process of intensity $n$ and the binomial point process.

Let $X_i$, $i \geq 1$, be a sequence of i.i.d. uniform random variables in $K(\en2)$
and put $\X_n:= \{X_1,\cdots,X_n\}$. For sake of simplicity, we denote by $f_k(\X_n \cap K(\en2))$ the number of $k$-dimensional faces of the convex hull of $\X_n$. In particular, we have
$$f_k(\X_n \cap K(\en2)):=
\sum_{X_i \in \X_n \cap K(\en2)} \xi_k(X_i, \X_n \cap K(\en2)).$$
We
start with two preliminary lemmas which describe  the growth of
$f_k(\X_n\cap K(\en2))$.

\begin{lemm} \label{Lem1} For all $k \in \{0,1,...,d-1 \}$ there is
a set $F(n), P[F(n)^c] = O(n^{-4d})$, and a constant $C_1 \in
(0, \infty)$ such that on $F(n)$ \be \label{JY4} |f_k(\X_n \cap
K(\en2)) - f_k(\X_{n + 1} \cap K(\en2))| \leq C_1 (\log n)^{k + 1}.
\ee
\end{lemm}

\noindent{\em Proof}.  
As in the proof of Lemma \ref{L4.0} and as on the pages 499-502 of \cite{Re}, there is a set $F_1(n)$ with
$P[F_1(n)^c] = O(n^{-4d}),$ such that on $F_1(n)$ we have for
$X_i \in K(\en2), \ 1 \leq i \leq n + 1,$
$$
\xi_{k}(X_i, \X_n \cap K(\en2)) = \xi_k(X_i, \X_n \cap K(\en2)
\cap B_{D_1 \epsilon_n}(X_i) ).$$ It follows that if $X_i \in
B^c_{D_1 \epsilon_n}(X_{n +1})  \cap K(\en2)$, then on $F_1(n)$ we
have
$$
\xi_k(X_i, \X_n \cap K(\en2)) = \xi_k(X_i, \X_{n + 1} \cap K(\en2)).
$$
Thus on $F_1(n)$ we have
$$
|f_k(\X_n \cap K(\en2)) - f_k(\X_{n + 1} \cap K(\en2))| $$
$$
\leq \xi_k(X_{n + 1}, \X_{n + 1}) +  \sum_{X_i \in B_{D_1
\epsilon_n}(X_{n + 1}) \cap K(\en2)} | \xi_k(X_i, \X_n \cap
K(\en2)) - \xi_k(X_i, \X_{n + 1} \cap K(\en2))|.
$$

The Lebesgue measure of $B_{D_1 \epsilon_n}(X_{n + 1}) \cap K(\en2)$
is $O( \epsilon_n^{d-1} \en2) = O( \epsilon_n^{d + 1} ) = O( \log
n/n)$.  There is thus a set $F_2(n)$, with $P[F_2^c(n)] = O(n^{-4d})$, such that on $F_2(n)$ we have
$$
\text{card} \{
 \X_n \cap B_{D_1 \epsilon_n}(X_{n + 1}) \cap K(\en2) \} = O (\log n).
$$

The proof of Lemma \ref{mom} shows that for $X_i \in B_{D_1
\epsilon_n}(X_{n + 1}) \cap K(\en2)$ there is a set $F_3(n),
P[F_3(n)^c] = O(n^{-4d })$, such that on $F_3(n)$ we have
$$
\xi_k(X_i, \X_{n}) = O(( \log n)^k).
$$

The same occurs for $\xi_k(X_{n+1},\X_{n+1})$. Now
on the set $F(n): = F_1(n) \cap F_2(n) \cap F_3(n)$ we get
\eqref{JY4}, concluding the proof of Lemma \ref{Lem1}.  \qed

\begin{lemm} \label{Lem2}
For all $k \in \{0,1,...,d-1 \}$ there is a constant $C_2$ such that
for all integers $l = 1,2,...,n$ we have
$$
P[|f_k(\X_n \cap K(\en2)) - f_k(\X_{n + l} \cap K(\en2))| \geq C_2 l
(\log n)^{k + 1}] \leq C_2 l n^{-4d}.$$
\end{lemm}

\noindent{\em Proof.}  We have
$$
|f_k(\X_n \cap K(\en2)) - f_k(\X_{n + l} \cap K(\en2))| \leq \sum_{i
= 0}^{l - 1} |f_k(\X_{n + i} \cap K(\en2)) - f_k(\X_{n + i + l} \cap
K(\en2))|.$$ By Lemma \ref{Lem1}, the $i$th summand is bounded by
$C_1 (\log (n + i))^{k + 1}$ on a set whose complement probability
is $O(n^{-4d})$.  Since $C_1 (\log (n + i))^{k + 1} \leq C (\log
2n)^{k + 1}$, the result follows.  \qed

\vskip.5cm

For every $\lambda>0$, let $N(\la)$ denote a Poisson variable of mean $\la$ and for
 every integer $n$ and $p\in (0,1)$, let $\rm{Bi}(n,p)$ denote a Binomial variable
 of parameters $n$ and $p$. 
The next result yields Theorem \ref{Th1b}.
\begin{theo} \label{dePo} For all $k \in \{0,1,...,d-1 \}$ we have
$$
|\Var f_k(K_n') - \Var f_k(K_{N(n)}')| = O\left( n^{1 - \frac{3}{d +
1} + o(1)} \right).
$$
\end{theo}
{\em Proof.}
For all integers $m$
we put $H_m := f_k(K_m')$. We have
$$
\Var H_n = \Var H_{N(n)} + \Var(H_n - H_{N(n)}) + 2 \Cov(H_{N(n)}, H_n - H_{N(n)} ).
$$
By \eqref{main1}, we have
$$
\Cov(H_{N(n)}, H_n - H_{N(n)} )  $$ $$\leq \sqrt{   \Var H_{N(n)}  }
 \  ||H_n - H_{N(n)} ||_2 = O \left( n^{(d-1)/2(d + 1)}  \right) ||H_n - H_{N(n)}
  ||_2
$$
It is thus enough to show \be \label{JY2}
  || H_n - H_{N(n)} ||_2^2 = O(n^{1 - \frac{4}{d + 1} + o(1)} ).
 \ee

Given the binomial and Poisson distributions ${\cal L}(\rm{Bi}(n, \en2))$ and ${\cal L}(N(n \en2))$, there exist
coupled random variables $\rm{Bi}(n, \en2)$ and $N(n \en2)$ such that
\be
\label{JY3} P[ \rm{Bi}(n, \en2) \neq N(n \en2)] \leq \en2;\ee
see e.g. (1.4) and (1.23) of \cite{BHJ}.

Enumerate the points $\P_n \cap K(\en2)$ by $X_1, X_2,...,X_{N(\en2)}$.
Given $\rm{Bi}(n, \en2)$, consider the coupled point set $\Y_n$ obtained by discarding or adding
i.i.d. points $X_i$ in $K(\en2)$:

 \be
 {\cal Y}_n : = \begin{cases} X_1,...,X_{N(\en2)- (N(\en2) - \rm{Bi}(n, \en2) )^+}
 , & {\rm if} \ N(\en2) \geq \rm{Bi}(n, \en2)
\\
X_1,...,X_{N(\en2) + (\rm{Bi}(n, \en2) - N(\en2))^+}
 , & {\rm if} \ N(\en2) < \rm{Bi}(n, \en2).
\end{cases}
 \ee
Then ${\cal Y}_n \eqd \X_n \cap K(\en2) = X_1, X_2,...,X_{\rm{Bi}(n, \en2)}.$  We use this coupling of the point sets
$\P_n \cap K(\en2)$ and $\X_n \cap K(\en2)$ in all that follows.

Denoting the convex hull of $m$ i.i.d. points $X_1,...,X_m$ on
$K(s)$ by $K(s)'_m$, we have
$$
|| H_n - H_{N(n)} ||_2^2 = \int (f_k(K_n') - f_k(K_{N(n)}'))^2 dP $$
$$= \int \left[ f_k(K(\en2)_{ {\rm{Bi}}(n, {\rm{vol}}(K(\en2))}') - f_k (K(\en2)_{N(n {\rm{vol}}(K(\en2)))}')
\right]^2 dP  + o(1),$$
where the last equality follows from the
$O(n^{-4d})$ probability bounds of Lemma 4.1(a), the bounds
$f_k(K(\en2)_{ j}') \leq C_3j^{d/2}$, as well as a standard
application of the Cauchy-Schwarz inequality. Let $E_n :=
\{\rm{Bi}(n,\rm{vol}(K(\en2))) \neq N(n \rm{vol}(K(\en2))) \}$ and recall from 
\eqref{JY3} that $P[E_n] \leq \en2$. On $E_n^c$ the integrand
vanishes. Thus
$$
|| H_n - H_{N(n)} ||_2^2 = \int \left[ f_k(K(\en2)_{ \rm{Bi}(n,\rm{vol}(K(\en2)))} ) - f_k (K(\en2)_{N(n \rm{vol}(K(\en2)))} ) \right]^2 {\bf 1}(E_n) dP +
o(1).
$$
By the Bernstein  inequality there is a constant $C_4$ such that for
all $p \in (0, 1/2)$ we have
$$
|{\rm{Bi}}(n, p) - np| \leq  C_4 (\log (np)) \sqrt{np} $$
with probability at least $1 - O(n^{-4d})$.  By Proposition A.2.3 of
\cite{BHJ}, and taking $C_4$ larger if necessary, we also have
$$
|{\rm{N}}(np) - np| \leq C_4 (\log (np)) \sqrt{np}$$ with
probability at least $1 - O(n^{-4d})$. A modification of Lemma
\ref{Lem2} shows that there is a set $G_n(1)$ with probability at
least $1 - O((\log n)^{1 +1/(d+1)}n^{1/2-1/(d+1)-4d})$ such that on
$G_n(1)$ we have
$$
|f_k(K(\en2)_{ \rm{Bi}(n,\rm{vol}(K(\en2)))} ) - f_k (K(\en2)_{n \lfloor\rm{vol}(K(\en2))\rfloor} )|^2 = O(
(\log n)^{2k + 4}  n \en2).
$$
Similarly, there is a set $G_n(2)$ with probability at least $1 -
O((\log n)^{1 +1/(d+1)}n^{1/2-1/(d+1)-4d})$ such that on $G_n(2)$ we
have
$$
|f_k(K(\en2)_{ \rm{N}(n \rm{vol}(K(\en2)))} ) - f_k (K(\en2)_{n \lfloor\rm{vol}(K(\en2))\rfloor} )|^2 = O(
(\log n)^{2k + 4}  n \en2).
$$
On the set  $G_n:= G_n(1) \cup G_n(2)$ we have \be \label{JY5}
|f_k(K(\en2)_{ \rm{Bi}(n, \en2)} ) - f_k (K(\en2)_{ \rm{N}(n \en2)}
)|^2 = O( (\log n)^{2k + 4}  n \en2). \ee

By McClullen's bound \cite{Mc}
$$
|f_k(K(\en2)_{ \rm{Bi}(n, \en2)} ) - f_k (K(\en2)_{\rm{N}(n \en2)}
)|^2 \leq C_3 ( \rm{Bi}(n, \en2)^{d} + \rm{N}(n \en2)^{d})$$ always
holds. It follows by the Cauchy-Schwarz inequality that $$ \int
\left[ f_k(K(\en2)_{ \rm{Bi}(n, \en2)} ) - f_k (K(\en2)_{N(n \en2)}
) \right]^2 {\bf 1}(E_n) {\bf 1}(G^c_n)dP = o(1),$$ whence in view
of \eqref{JY5}
$$
|| H_n - H_{N(n)} ||_2^2 = O \left(  (\log n)^{2k + 4} n \en2 \int
{\bf 1}(E_n) {\bf 1}(G_n)dP \right) + o(1).$$

It follows that
$$
|| H_n - H_{N(n)} ||_2^2 = O((\log n)^{2k + 4} n \en2 P[E_n]) + o(1)
= O((\log n)^{2k + 4} n \epsilon_n^4 ) + o(1).$$

This shows \eqref{JY2} and concludes the proof of Theorem
\ref{dePo}. \qed

Pierre Calka, Laboratoire de Math\'ematiques Rapha\"el Salem, Universit\'e de Rouen,
 Avenue de l'Universit\'e, BP.12, Technop\^ole du Madrillet, F76801 Saint-Etienne-du-Rouvray France;
 \\   \ \ \  {\texttt pierre.calka@univ-rouen.fr}

J. E. Yukich, Department of Mathematics, Lehigh University,
Bethlehem PA 18015;\\
 \ \ {\texttt joseph.yukich@lehigh.edu}

\end{document}